\documentclass[10pt]{amsart}
 
\usepackage{amsmath,amsthm}
\usepackage{amssymb}
\usepackage{euscript}

\input epsf

\newcommand{\sqsp}{\renewcommand{\baselinestretch}{1.3}\tiny\normalsize}

\numberwithin{equation}{section}  % [if desired]
 
\newtheorem{thm}[equation]{Theorem}
\newtheorem{prop}[equation]{Proposition}
\newtheorem{cor}[equation]{Corollary}
\newtheorem{lemma}[equation]{Lemma}

\theoremstyle{definition}
\newtheorem{defn}[equation]{Definition}

\newtheorem{rmk}[equation]{Remark}

\begin{document}

\title[Khovanov's invariants for alternating knots]{The support of the Khovanov's invariants for alternating knots}
\author{Eun Soo Lee}

\begin{abstract}
In this article, we prove the conjecture of Bar-Natan, Garoufalidis, and Khovanov's on the support of the Khovanov's invariants for alternating knots.
\end{abstract}

\date{\today}
% \subjclass{57M27}
% \keywords{}

\address{Department of Mathematics, MIT, Cambridge, MA 02139, U.S.A.}
\email{eslee@math.mit.edu}

\maketitle
 
%%%%%%%%%%%%%%%%%%%%%%%%%%

\sqsp

\bigskip
\section{Introduction}

   In \cite{K}, Khovanov constructed invariants of knots and links, and then, in \cite{B} and \cite{G}, Bar-Natan, Garoufalidis, and Khovanov formulated conjectures on the values of Khovanov's invariants for alternating knots.
   This article provides a proof of one of those conjectures.

   The following is the main theorem.

\begin{thm}
\label{thm:main}
   For any oriented non-split alternating link $L$, $Kh(L)$ is supported in two lines $\deg(q) = 2 \deg(t) - \sigma(L) \pm 1 $, its nonzero coefficient of the smallest degree in $t$ is on the line $\deg(q) = 2 \deg(t) - \sigma(L) - 1 $, its nonzero coefficient of the largest degree in $t$ is on the line $\deg(q) = 2 \deg(t) - \sigma(L) + 1 $, and those coefficients are 1.
   In other words,
\[ Kh(L) = \sum_{i=p}^m ( a_i t^i q^{2i - \sigma(L) - 1} + b_i t^i q^{2i - \sigma(L) + 1} ) \]
for some $p \leq m$ with $a_p = b_m = 1$.
\end{thm}

\begin{rmk}
   In the formula above, $\deg(q) = 2 \deg(t) - \sigma(L) - 1 $ line, and $\deg(q) = 2 \deg(t) - \sigma(L) + 1 $ line will be called the diagonal, and the subdiagonal, respectively, and the positions of $a_p = 1 $ and $b_m = 1$ will be referred to as the top at $(p, 2p - \sigma(L) - 1)$ and the bottom at $(m, 2m - \sigma(L) + 1)$, thinking of the table of coefficients in which the powers of $t$ increase from left to right, and the powers of $q$ increase from top to bottom. These terms will be applied to Khovanov's cohomology groups as well.

   We follow the convention of defining $\sigma(L)$ to be the signature of the sum of the Seifert matrix of $L$ and its transpose, which is the negative of the signature in \cite{B}.

   For the sign of a crossing, we follow the convention indicated below, which is opposite to what is used in \cite{K}, so that $x(D)$ and $y(D)$ for an oriented diagram $D$ are the number of negative crossings and positive crossings, respectively.

   We follow the notations and terminologies defined in \cite{K}. Consult \cite{K} for undefined terminologies and notations.
\end{rmk}

\[
\epsfysize=0.75in\epsfbox{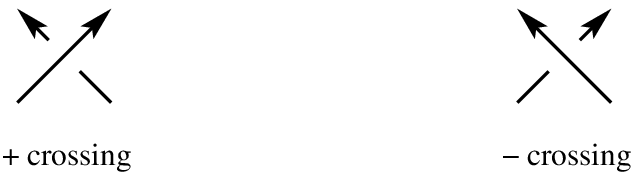}
\]

\bigskip
\section{Properties of black and white coloring of an alternating link diagram}

   Let $D$ be a link diagram, which is a regular projection of a link together with the information of relative height at each double point. 
For brevity of the statements to follow, let's think of the diagram on $S^2$ rather than $\mathbb{R}^2$. The regions of $S^2$ divided by $D$ can be colored black and white in checkerboard fashion. 

   At each crossing, a coloring of the nearby regions falls into one of the two following patterns. 

\[
\epsfysize=0.69in\epsfbox{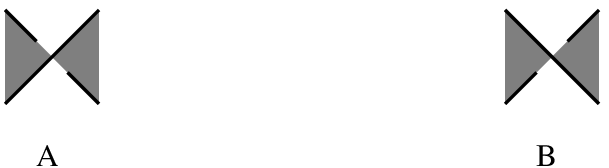}
\]  
   
   If $D$ is alternating, then adjacent crossings have the same coloring pattern.

\[
\epsfysize=0.69in\epsfbox{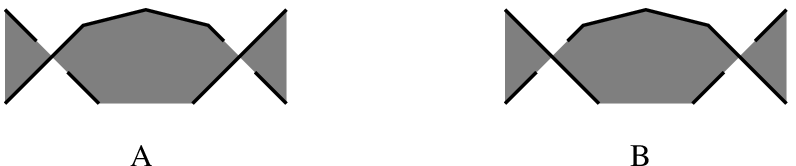}
\]  
   
   Thus, in a coloring of a non-split alternating diagram $D$, only one of the pattern A or B appears for every crossing. Reversion of a coloring changes that pattern. 

   To deal with the resolutions of $D$, consider the resolutions of a colored diagram as below.
   
\[
\epsfysize=1.50in\epsfbox{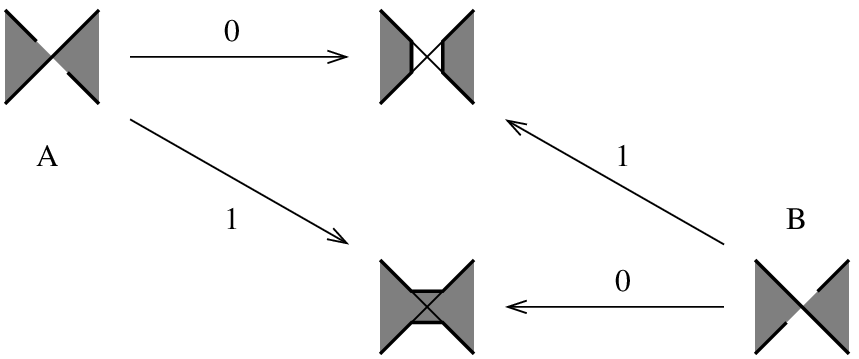}
\]

\begin{defn}
   For a non-split alternating diagram $D$, \emph{the coloring} of $D$ is the coloring of $D$ in which only pattern A appears. \emph{The coloring} of a resolution of $D$ is the coloring of that resolved diagram induced from the coloring of $D$.   
\end{defn} 

   For the coloring of $D(\emptyset)$ (0-resolutions of pattern A), the trace of each crossing lies in a white region. Now, our claim is:

\begin{prop}
\label{prop:goodD}
   For a reduced non-split alternating diagram $D$, the components of $D(\emptyset)$ bound non-overlapping black disks in the coloring of it. Each black disk corresponds to each of the black regions in the coloring of $D$. Furthermore, every pair of black disks are connected by a chain of black disks, which are connected by the trace of the crossings of $D$. Also, no trace of crossing connects a black disk to itself. 
\end{prop}

   Here is a visualization of our claim for the left-handed trefoil and a figure 8 knot. (The unbounded black region shown below is a disk in $S^2$.) 

\[
\epsfysize=0.94in\epsfbox{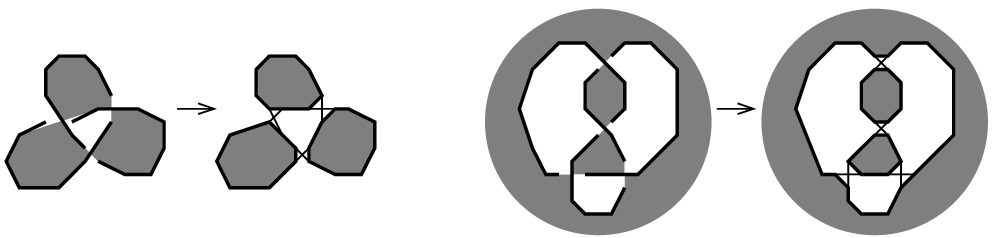}
\]  

\begin{proof}
   At each crossing, the 0-resolution separates incident black regions. So, we get  correspondence between the black regions in the coloring of $D$ and those in the coloring of $D(\emptyset)$. (While the white regions of $D$ merge in the process.) 
    
   In the coloring of $D(\emptyset)$, there's no trace of crossing in black regions. That implies:

   $\bullet$ if there is a black region which is not a disk, then $D$ is split. 

   $\bullet$ if there is a pair of black disks which cannot be connected by any chain, then $D$ is split. 

   $\bullet$ if there is a trace of crossing connecting a black disk to itself, that crossing is removable, so $D$ is not reduced.  

\[
\epsfysize=1.00in\epsfbox{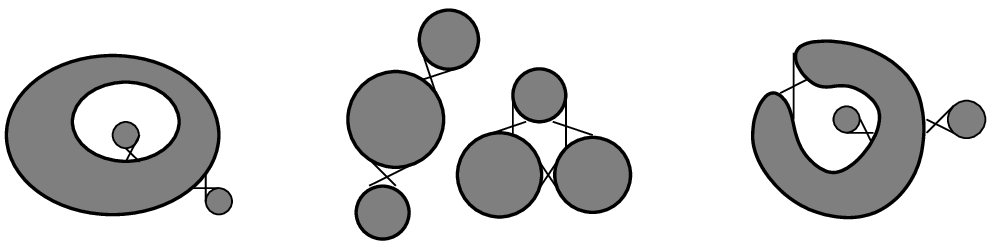}
\]  

\end{proof}

\begin{defn}   
   For a link diagram $D$, let $c(D)$ be the number of crossings of $D$, and $o(D)$ be the number of components of $D(\emptyset)$.
\end{defn}

   For $D$ in proposition \ref{prop:goodD}, $o(D)$ agrees with the number of black disks in the coloring of $D(\emptyset)$.

   Let $\mathcal{I}$ be an ordered set of crossings of $D$. Note that $D(\mathcal{I})$ agrees with $D^!(\emptyset)$, and that 
$ o(D) + o(D^!) $ equals the total number of black and white regions in the coloring of $D$, which is $ c(D) + 2 $.

   We need one further step for the inductive argument to be used in our proof of the main theorem.

\begin{prop}
\label{prop:betterD}
   Let $D$ be a reduced non-split alternating link diagram with $c(D) > 0$.
Then one of the following holds.  

   \textsl{I}. There is a pair of black disks in the coloring of $D(\emptyset)$ connected by exactly one crossing.

   \textsl{II}. There is a pair of black disks in the coloring of  $D^!(\emptyset)$ connected by exactly one crossing.

   \textsl{III}. $D$ is a connected sum of $D'$ and the simplest link, for another reduced non-split alternating link diagram $D'$ with $c(D) - 2$ crossings.
 
\[
\epsfysize=0.75in\epsfbox{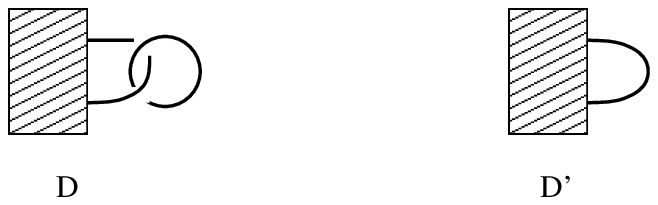}
\]  

\end{prop}

\begin{proof}
   Since $o(D) + o(D^!) = c(D) + 2$, one of the following holds.
   
   \textsl{i}. $o(D) > c(D)/2 + 1$.
   
   \textsl{ii}. $o(D^!) > c(D)/2 + 1$.
   
   \textsl{iii}. $ o(D) = o(D^!) = c(D)/2 + 1$.

[\textsl{i} $\Rightarrow$ \textsl{I}]
   For $o(D)$ black disks to be connected to each other by chains of connected disks, there are at least $o(D) -1$ different pairs that are connected by crossings. If $2 (o(D) -1) > c(D) $, then at least one of those pairs is connected by exactly one crossing.   
   
[\textsl{ii} $\Rightarrow$ \textsl{II}]
   Same as \textsl{i} $\Rightarrow$ \textsl{I}.

[\textsl{iii} \& not \textsl{I} \& not \textsl{II} $\Rightarrow$ \textsl{III}]
   To fail \textsl{I}, there are exactly $o(D) - 1$ different pairs that are connected by crossings and those pairs are connected by exactly two crossings. 

   Consider a graph consists of $o(D)$ vertices and $o(D) -1$ edges. Each vertex represents each black disk. For each pair of black disks connected by two crossings, there's an edge joining the corresponding pair of vertices.   
This graph is connected, so it is a tree.

   For an edge $\{a,b\}$, mark the $a$-end of it with arrow if the two crossings connecting the disk $a$ and $b$ is adjacent on the boundary of $a$.
For example,

\[
\epsfysize=1.10in\epsfbox{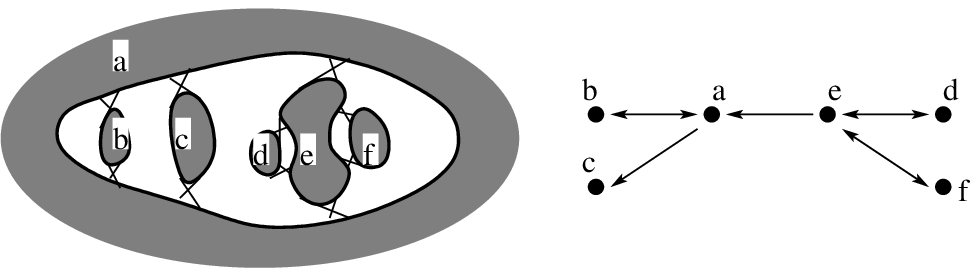}
\]  

   A vertex of a tree is called a pendent vertex if it is incident with only one edge, and an edge is called a pendent edge if it is incident with a pendent vertex. 
   If $a$ is a pendent vertex, the unique edge incident with $a$ is necessarily marked at the $a$-end. If $b$ is not a pendent vertex, at least two edges have marked $b$-end, since the two crossings connecting disks $b$ and $c$ and those connecting disks $b$ and $d$ never alternate.   
     
\[
\epsfysize=0.69in\epsfbox{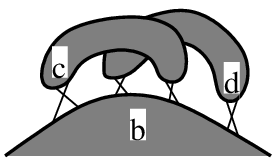}
\]  

   If $o(D)=2$, then there's only one edge, that is a pendent edge, and both ends of that edge is marked. For $o(D) > 2$, let $p$ be the number of the pendent vertices. The number of the pendent edges is also $p$. There are at least $p + 2(o(D)-p)$ marked ends, but the number of non-pendent edges is $o(D)- 1- p$, so there is at least one pendent edge with both ends marked.
That implies \textsl{III} (up to relocation of $\infty$).

\[
\epsfysize=0.81in\epsfbox{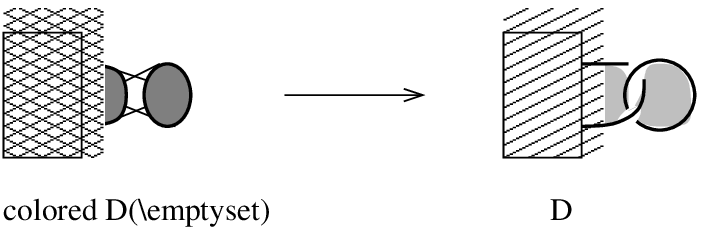}
\]  
   
\end{proof}

\begin{lemma}
\label{lemma:box-support}
   For a reduced non-split alternating diagram $D$, $\overline{\mathcal{H}}^{i,j}(D)$ is supported in the box $0 \leq i \leq c(D)$ and $-o(D) \leq j \leq 2c(D)-o(D)+2$,
with $\overline{\mathcal{H}}^{0,-o(D)}(D) = \overline{\mathcal{H}}^{c(D),2c(D)-o(D)+2}(D) = \mathbb{Z}$.
\end{lemma}

\begin{proof}
   First of all, it is clear from the construction of $\overline{\mathcal{C}}(D)$ that $\overline{\mathcal{C}}^{i,j}(D) = 0$ unless $0 \leq i \leq c(D)$
 
   When a resolution of $D$ is changed to another resolution of $D$ by replacing one 0-resolution by 1-resolution, the number of components either increases by 1 or decreases by 1. That insures $\overline{\mathcal{C}}^{i,j}(D)$ to be supported in $-o(D) \leq j \leq 2c(D)-o(D)+2$.

   Proposition \ref{prop:goodD} implies that $D(\emptyset)$ has one more component than any $D(a)$ has, because two black disks merge into one in the process. In terms of $\overline{\mathcal{C}}^{i,j}(D)$, this means 
\[\overline{\mathcal{C}}^{i,j}(D) = 
\begin{cases}
\mathbb{Z} & \textrm{ if } i=0, j = -o(D) \\
0 & \textrm{ if } i > 0, j = -o(D)\textrm{ ,}
\end{cases}\]        
so the half of the result follows.

   For the other half, look at the other end. $D(\mathcal{I})= D^!(\emptyset)$ also has one more component than any $D(\mathcal{I}-\{a\})= D^!(a)$ has, so that
\[\overline{\mathcal{C}}^{i,j}(D) = 
\begin{cases}
\mathbb{Z} & \textrm{ if } i=c(D), j = 2c(D) -o(D)+2 \\
0 & \textrm{ if } i < c(D), j = 2c(D)-o(D)+2 \textrm{ .}
\end{cases}\]      
\end{proof}

   Let $D$ be a diagram satisfying \textsl{I} in proposition \ref{prop:betterD}. Let $a$ be a crossing of $D$ connecting a pair of black disks that no other crossing connects. Order $\mathcal{I}$ so that $a$ comes the last. Let $D(*0)$, $D(*1)$ be the diagram obtained from $D$ by resolving $a$ to its 0- and 1-resolution, respectively. $D(*0)$ still has the property that $D(*0)(\emptyset)$ has one more component than any $D(*0)(b)$ has. The use of \textsl{I} is that it allows $D(*1)$ to have that property, too. 

\begin{cor}
\label{cor:box-support}
   In the above setting, $\overline{\mathcal{H}}^{i,j}(D(*0))$ is supported in the box $0 \leq i \leq c(D(*0))$ and $-o(D(*0)) \leq j \leq 2c(D(*0))-o(D(*0))+2$,
with $\overline{\mathcal{H}}^{0,-o(D(*0))}(D(*0)) = \mathbb{Z}$, and 
$\overline{\mathcal{H}}^{i,j}(D(*1))$ is supported in the box $0 \leq i \leq c(D(*1))$ and $-o(D(*1)) \leq j \leq 2c(D(*1))-o(D(*1))+2$,
with $\overline{\mathcal{H}}^{0,-o(D(*1))}(D(*1)) = \overline{\mathcal{H}}^{c(D(*1)),2c(D(*1))-o(D(*1))+2}(D(*1)) = \mathbb{Z}$. 
\end{cor}
   
   Finally, to apply induction hypothesis to $D(*0)$ and $D(*1)$ later on, they need to be non-split alternating.

\begin{prop}
   In the above setting, $D(*0)$ and $D(*1)$ are non-split alternating.
\end{prop}

\begin{proof}
   Alternating property is easy to see.
   
   To be non-split, their black disks in the induced coloring have to be connected. That is clear for $D(*1)$. For $D(*0)$, if the black disks of $D(\emptyset)$ are disconnected after removing $a$, then $a$ were a removable crossing in $D$, contradiction to $D$ being reduced.
\end{proof}

\bigskip
\section{Signature of an alternating link}

   This section consists of the result of \cite{GL} and an application to alternating links, to relate the place of diagonal with the signature in the main theorem.
 
\begin{defn}[Goeritz matrix : following \S1 of \cite{GL}]
   Let $D$ be an oriented link diagram. Color the regions of $\mathbb{R}^2$ (or $S^2$) divided by $D$ in checkerboard fashion. Denote the white regions by $X_0, X_1, \cdots, X_n$. Assume that each crossing is incident to two distinct white regions. Assign an incidence number $\eta(a)= \pm 1$ to each crossing $a$ as in the figure below. For $0 \leq i , j \leq n $ define 
\[g_{ij} = 
\begin{cases} 
- \sum_{a \textrm{ incident to both }X_i \textrm{ and }X_j} \eta(a)  & \textrm{ for }i \neq j \\
- \sum_{0 \leq k \leq n , k \neq i} g_{ik}  & \textrm{ for }i = j \textrm{ .}
\end{cases}\]

   The Goeritz matrix $G(D)$ of $D$ is the $n \times n $ (not $ (n+1) \times (n+1)$ !) symmetric matrix $G(D)=(g_{ij})_{1 \leq i ,j \leq n}$. 
\end{defn}

\[
\epsfysize=0.75in\epsfbox{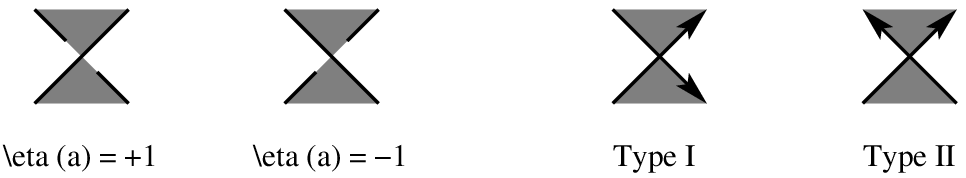}
\]  

   The signature of an oriented link can be obtained from the signature of Goeritz matrix of its diagram by adding a correction term.

\begin{thm}[Theorem 6 in \cite{GL}]
   For an oriented link $L$, 
\[ \sigma(L) = \mathrm{sign}\ G(D) - \mu(D) \]
for its diagram $D$, where 
$ \mu(D) = \sum \eta(a) $, summed over all crossings $a$  of type II. (See the figure above)
\end{thm}

\begin{prop}
\label{prop:shift}
   For an oriented non-split alternating link $L$ and a reduced alternating diagram $D$ of $L$, $\sigma(L) = o(D) - y(D) -1$.
\end{prop}

\begin{proof}
    The non-split alternating property of $D$ implies either $\eta(a)=1$ for all crossings $a$ or $\eta(a)= -1$ for all $a$. By reversing the coloring if necessary, we may assume that $\eta(a)=1$ for all $a$. 

   Since $D$ is non-split alternating, components of the resolution $D(\emptyset)$ bound non-overlapping disks. Our arrangement $\eta(a)=1$ forces those disks to be correspondent with white regions, $+$ crossings to be of the type II, and $-$ crossings to be of the type I. Reducedness of $D$ insures that each crossing is incident to two distinct white regions.

   Then, $G(D)$ is a positive definite matrix of rank $o(D)-1$ and $\mu(D) = y(D)$, so that $\sigma(L) = \mathrm{sign}\ G(D) - \mu(D) = o(D) -1 - y(D)$.  
\end{proof}

\bigskip
\section{Proof of the main theorem}   

   The proof is based on induction on the number of crossings of a link diagram. First, we prove the theorem for some number $s(L)$ instead of $\sigma(L)$, and then, show $s(L)=\sigma(L)$.

\begin{lemma}
\label{lemma:H->H}
   The chain complexes $\overline{\mathcal{C}}(D)$, $\overline{\mathcal{C}}(D(*0))$, and $\overline{\mathcal{C}}(D(*1))[-1]\{-1\}$ form a short exact sequence
\[ 0 \rightarrow \overline{\mathcal{C}}(D(*1))[-1]\{-1\} \rightarrow
\overline{\mathcal{C}}(D) \rightarrow \overline{\mathcal{C}}(D(*0)) \rightarrow 0 \]
with degree preserving maps, 
so that $\overline{\mathcal{H}}(D)$ is an extension of the kernel and cokernel of the connecting map $\delta$ as a bigraded $\mathbb{Z}$-module. 
   \[ 0 \rightarrow \mathrm{Coker} \ \delta
\rightarrow \overline{\mathcal{H}}(D) \rightarrow \mathrm{Ker} \ \delta \rightarrow 0 \] 

   In particular, the support of $\overline{\mathcal{H}}(D)$ is included in the union of the support of $\overline{\mathcal{H}}(D(*0))$ and $\overline{\mathcal{H}}(D(*1))[-1]\{-1\}$.
\end{lemma}

\begin{proof}
   Consider the $n$-dimensional cube associated with $\overline{\mathcal{C}}(D)$. The $(n-1)$-dimensional sub-cube indexed by the subsets of $\mathcal{I}$ not containing $a$ is associated with $\overline{\mathcal{C}}(D(*0))$, and the $(n-1)$-dimensional sub-cube indexed by the subsets of $\mathcal{I}$ containing $a$ is associated with  $\overline{\mathcal{C}}(D(*1))[-1]\{-1\}$. 
   Thus, $\overline{\mathcal{C}}(D)$ decomposes into 
\[\overline{\mathcal{C}}(D) = \overline{\mathcal{C}}(D(*0)) \oplus \overline{\mathcal{C}}(D(*1))[-1]\{-1\} \textrm{ .}\]

   Let $d$, $d_0$, and $d_1$ be the coboundary maps of the complexes $\overline{\mathcal{C}}(D)$, $\overline{\mathcal{C}}(D(*0))$, and $\overline{\mathcal{C}}(D(*1))[-1]\{-1\}$, respectively. 
   Since no index for the sub-cube associated with $\overline{\mathcal{C}}(D(*1))[-1]\{-1\}$ is contained in an index for the sub-cube associated with $\overline{\mathcal{C}}(D(*0))$, 
$d$ doesn't have any part from  $\overline{\mathcal{C}}(D(*1))[-1]\{-1\}$ to  $\overline{\mathcal{C}}(D(*0))$, and hence, $d$ decomposes into 
\[ d(z,w) = (d_0 z, \xi z - d_1 w) \textrm{ ,}\]
where $\xi : \overline{\mathcal{C}}(D(*0)) \to \overline{\mathcal{C}}(D(*1))[-1]\{-1\}$ is the part of the coboundary map $d$ from $\overline{\mathcal{C}}(D(*0))$ to $ \overline{\mathcal{C}}(D(*1))[-1]\{-1\}$.

   Now, it is easy to see that
\[ 0 \rightarrow \overline{\mathcal{C}}(D(*1))[-1]\{-1\} \rightarrow
\overline{\mathcal{C}}(D) \rightarrow \overline{\mathcal{C}}(D(*0)) \rightarrow 0 \]
is a short exact sequence of chain complexes (after a little adjustment of sign), and that  $\delta : \overline{\mathcal{H}}(D(*0)) \to \overline{\mathcal{H}}(D(*1))[-1]\{-1\}$ is induced by $\xi$.
\end{proof}

\begin{thm}
\label{thm:almost}
   For any non-split alternating link diagram $D$, $\overline{\mathcal{H}}^{i,j}(D) \otimes_{\mathbb{Z}} {\mathbb{Q}}$ is supported in two lines $ j = 2i - s \pm 1 $ for some integer $s$ with the top and bottom on the diagonal and the subdiagonal, respectively.
\end{thm}

\begin{proof}
   For the base case, the theorem holds for the unknotted diagram of unknot.
      
   Assume that the statement is true for all such diagrams with less than $n$ crossings. Let $D$ be a non-split alternating link diagram with $n$ crossings. If $D$ is not reduced, then $\overline{\mathcal{H}}(D)$ is a shift of $\overline{\mathcal{H}}(D')$ for some such diagram $D'$ with less than $n$ crossings, so the statement is true for $D$ as well. 

   Let $D$ be reduced. By the duality theorem for mirror image (proposition 32 and corollary 11 in \cite{K}), it is enough to show for either $D$ or $D^!$. So, we may assume that $D$ has the property \textsl{I} or \textsl{III} in proposition \ref{prop:betterD}.
   
[\textsl{Case I}] 
   The induction hypothesis applies to $D(*0)$ and $D(*1)$.  $\overline{\mathcal{H}}(D(*0)) \otimes_{\mathbb{Z}} {\mathbb{Q}}$ is supported in two lines with the top at $(0, -o(D(*0)))$, and $\overline{\mathcal{H}}(D(*1)) \otimes_{\mathbb{Z}} {\mathbb{Q}}$ 
is also supported in two lines with the top at $(0, -o(D(*1)))$.

   Since $o(D)=o(D(*0))=o(D(*1))+1$, the diagonal and the subdiagonal of $\overline{\mathcal{H}}(D(*0)) \otimes_{\mathbb{Z}} {\mathbb{Q}}$ agree with those of $\overline{\mathcal{H}}(D(*1))[-1]\{-1\} \otimes_{\mathbb{Z}} {\mathbb{Q}}$.
   By lemma \ref{lemma:H->H} and lemma \ref{lemma:box-support}, $\overline{\mathcal{H}}(D) \otimes_{\mathbb{Z}} {\mathbb{Q}}$ is supported in two lines with the top at $(0, -o(D))$ and the bottom at $(n, 2n-o(D)+2)$.

[\textsl{Case III}]
   Our $D(*0)$ and $D(*1)$ are as below, and the induction hypothesis applies to $D'$.
   
\[
\epsfysize=0.88in\epsfbox{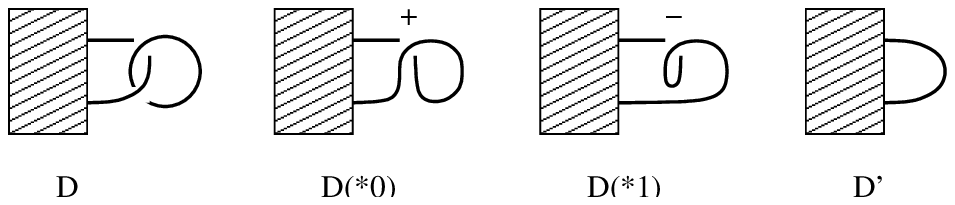}
\]  
   
   Orient $D'$, $D(*0)$ and $D(*1)$ accordingly. $\overline{\mathcal{H}}(D')$, $\overline{\mathcal{H}}(D(*0))$ and $\overline{\mathcal{H}}(D(*1))[-1]\{-1\}$ are shift of each other as follows.
\begin{eqnarray*}
\overline{\mathcal{H}}(D(*0))
&=& \mathcal{H}(D(*0))[-x(D(*0))]\{-2x(D(*0)) + y(D(*0))\}  \\
&=& \mathcal{H}(D')[-x(D')]\{-2x(D') + y(D')+ 1\} \\
&=& \overline{\mathcal{H}}(D')[0]\{1\}
\end{eqnarray*}
\begin{eqnarray*}
\overline{\mathcal{H}}(D(*1))[-1]\{-1\}
&=& \mathcal{H}(D(*1))[-x(D(*1))-1]\{-2x(D(*1)) + y(D(*1)) -1\}  \\
&=& \mathcal{H}(D')[-x(D')-2]\{-2x(D') + y(D') -3\} \\
&=& \overline{\mathcal{H}}(D')[-2]\{-3\}
\end{eqnarray*}

   By induction hypothesis, $\overline{\mathcal{H}}(D(*0)) \otimes_{\mathbb{Z}} {\mathbb{Q}}$ is supported in two lines with the top at $(0, -o(D')-1)$, and $\overline{\mathcal{H}}(D(*1))[-1]\{-1\} \otimes_{\mathbb{Z}} {\mathbb{Q}}$ 
is also supported in two lines with the top at $(2, -o(D')+ 3)$. Their diagonals and subdiagonals agree. 

   Again, by lemma \ref{lemma:H->H} and lemma \ref{lemma:box-support}, $\overline{\mathcal{H}}(D) \otimes_{\mathbb{Z}} {\mathbb{Q}}$ is supported in two lines with the top at $(0, -o(D)) = (0, -o(D')-1)$ and the bottom at $(n, 2n-o(D)+2)$.
\end{proof}

   Let $L$ be an oriented non-split alternating link and $D$ be a reduced alternating diagram of $L$. From theorem \ref{thm:almost}, we can conclude that $\mathcal{H}(D)  \otimes_{\mathbb{Z}} {\mathbb{Q}} = \overline{\mathcal{H}}(D)[x(D)]\{2x(D)-y(D)\} \otimes_{\mathbb{Z}} {\mathbb{Q}} $ has the top at $(-x(D), -2x(D)+y(D)-o(D))$. Since the top is on the diagonal, our $s(L)$ equals $ o(D) - y(D) -1$.

   In proposition \ref{prop:shift}, we saw that $\sigma(L) = o(D) - y(D) -1$. That finishes the proof of the main theorem. 

   We can also tell something about the support of $\mathrm{Tor}(\mathcal{H}(L))$. The duality theorem for mirror image (proposition 32 and corollary 11 in \cite{K}) still allows us to work with either one or its mirror image. Thus, the following could have been included in the induction.
   
\begin{cor}
   For an oriented non-split alternating link $L$ with the top of $Kh(L)$ at $(p, 2p - \sigma(L) - 1)$ and the bottom at $(m, 2m - \sigma(L) +1)$,
$\mathrm{Tor}({\mathcal{H}}^{i,j}(L)) = 0 $ unless $p+1 \leq i \leq m$ and 
$j = 2i -\sigma(L) -1 $. 
\end{cor}

%%%%%%%%%%%%%%%%%%%%%%%%%

\end{document}